\documentclass[12pt,leqno]{article}
\usepackage{amsfonts}
\usepackage{amsmath, amsthm, amsfonts, amssymb, color}
\usepackage{mathrsfs}
\usepackage{color}
\newtheorem{thm}{Theorem}[section]

\newtheorem{prp}[thm]{Proposition}

\theoremstyle{definition}

\newcommand{\scr}[1]{\mathscr #1}
\definecolor{wco}{rgb}{0.5,0.2,0.3}

\numberwithin{equation}{section} \theoremstyle{remark}

\newcommand{\ua}{\uparrow}

\newcommand{\var}{\textup{Var}}

\title{{\bf Weak Poincar\'e Inequalities for Convergence Rate of Degenerate Diffusion Processes }\footnote{Supported in
 part by  NNSFC  (11431014,11626245,11626250).} }
\author{{\bf Martin Grothaus$^{b)}$ and Feng-Yu Wang$^{a),c)}$ }\\
\footnotesize{a)  Center for Applied Mathematics, Tianjin University, Tianjin 300072, China}\\
\footnotesize{b) Mathematics Department,    Kaiserslautern University,   67653 Kaiserslautern, Germany}\\
 \footnotesize{c) Department of Mathematics,
Swansea University, Singleton Park, SA2 8PP, UK}\\
\footnotesize{grothaus@mathematik.uni-kl.de; wangfy@bnu.edu.cn, F.-Y.Wang@swansea.ac.uk}
 }

\begin{document}
\allowdisplaybreaks
\def\R{\mathbb R}  \def\ff{\frac} \def\ss{\sqrt} \def\B{\mathbf
B}
\def\N{\mathbb N} \def\kk{\kappa} \def\m{{\bf m}}
\def\ee{\varepsilon}\def\ddd{D^*}
\def\dd{\delta} \def\DD{\Delta} \def\vv{\varepsilon} \def\rr{\rho}
\def\<{\langle} \def\>{\rangle} \def\GG{\Gamma} \def\gg{\gamma}
  \def\nn{\nabla} \def\pp{\partial} \def\E{\mathbb E}
\def\d{\text{\rm{d}}} \def\bb{\beta} \def\aa{\alpha} \def\D{\scr D}
  \def\si{\sigma} \def\ess{\text{\rm{ess}}}
\def\beg{\begin} \def\beq{\begin{equation}}  \def\F{\scr F}
\def\Ric{\text{\rm{Ric}}} \def\Hess{\text{\rm{Hess}}}
\def\e{\text{\rm{e}}} \def\ua{\underline a} \def\OO{\Omega}  \def\oo{\omega}
 \def\tt{\tilde} \def\Ric{\text{\rm{Ric}}}
\def\cut{\text{\rm{cut}}} \def\P{\mathbb P} \def\ifn{I_n(f^{\bigotimes n})}
\def\C{\scr C}      \def\aaa{\mathbf{r}}     \def\r{r}
\def\gap{\text{\rm{gap}}} \def\prr{\pi_{{\bf m},\varrho}}  \def\r{\mathbf r}
\def\Z{\mathbb Z} \def\vrr{\varrho} \def\ll{\lambda}
\def\L{\scr L}\def\Tt{\tt} \def\TT{\tt}\def\II{\mathbb I}
\def\i{{\rm in}}\def\Sect{{\rm Sect}}  \def\H{\mathbb H}
\def\M{\scr M}\def\Q{\mathbb Q} \def\texto{\text{o}} \def\LL{\Lambda}
\def\Rank{{\rm Rank}} \def\B{\scr B} \def\i{{\rm i}} \def\HR{\hat{\R}^d}
\def\to{\rightarrow}\def\l{\ell}\def\iint{\int}
\def\EE{\scr E}\def\Cut{{\rm Cut}}
\def\A{\scr A} \def\Lip{{\rm Lip}}
\def\BB{\scr B}\def\Ent{{\rm Ent}} \def\Var{{\rm Var}}

\maketitle

\begin{abstract} For a contraction $C_0$-semigroup on   a separable Hilbert space, the decay rate is   estimated by using the weak Poincar\'e inequalities for  the symmetric and anti-symmetric part of the generator. As applications, non-exponential convergence rate is characterized for a class of degenerate diffusion processes, so that the study of hypocoercivity  is extended. Concrete examples are presented.
\end{abstract} \noindent
 AMS subject Classification:\  60H10, 37A25, 47D07.   \\
\noindent
 Keywords: Degenerate diffusion semigroup,   hypocercivity, weak Poincar\'e inequality, convergence rate.
 \vskip 2cm

\section{Introduction}

Let $(E,\F, \mu)$ be a probability space and let  $(\EE,\D(\EE))$ be the quadratic  form associated with  a Markov semigroup $P_t$ on $L^2(\mu)$.  The weak Poincar\'e inequality
\beq\label{WPS} \var_{\mu}(f):= \mu(f^2)-\mu(f)^2\le \aa(r)\EE(f,f)+\aa(r) \|f\|_{\rm osc}^2,\ \ r>0, f\in \D(\EE)\end{equation} with rate function $\aa: (0,\infty)\to (0,\infty)$ was introduced in \cite{RW01} to describe the following convergence rate of $P_t$ to $\mu$:
$$\xi(t):= \sup_{\|f\|_{\rm osc}\le 1} \Var_\mu(P_t f),\ \ t>0.$$  Explicit correspondence between $\aa$ and $\xi$ has been presented in \cite{RW01}. In particular, the weak Poincar\'e inequality
\eqref{WPS} is always available for elliptic diffusion processes. However, it does not hold when the Dirichlet form is reducible.  A typical example is the stochastic Hamiltonian system on $\R^d\times \R^d$:
\beq\label{E1} \beg{cases} \d X_t= Y_t\d t,\\
\d Y_t= \ss 2 \,\d B_t -\big(\nn^{(1)} V(X_t)+Y_t\big)\d t,\end{cases}\end{equation}
where $B_t$ is the Brownian motion on $\R^d$, $\nn^{(1)}$ is the gradient operator in the first component $x\in\R^d$,  and $V\in C^2(\R^d)$   satisfies
\beq\label{VV} \|\nn^2 V\|\le M(1+|\nn V|)\end{equation} for some constant $M>0$ and $Z(V):=\int_{\R^d}\e^{-V(x)}\d x<\infty$.
In this case the invariant probability measure of the diffusion process is  $\mu=\mu_1\times\mu_2$, where
$\mu_1(\d x)= Z(V)^{-1} \e^{-V(x)}\d x$ and $\mu_2$ is the standard Gaussian measure on $\R^d$. Let $\nn^{(2)}$ be the gradient operator in the second component $y\in \R^d$.
Then the associated energy form satisfies
$\EE(f,f)= \mu(|\nn^{(2)}f|^2)$, and is   thus reducible.

On the other hand, according to C. Villani   \cite{Vil09},   if the Poincar\'e inequality
\beq\label{P} \Var_{\mu_1}(f):= \mu_1(f^2)- \mu_1(f)^2 \le c_1 \mu_1(|\nn f|^2),\ \ f\in C_b^1(\R^d)\end{equation} holds for some constant $c_1>0$, then the Markov semigroup $P_t$ associated with \eqref{E1} converges exponentially to $\mu$ in the sense that
$$\mu\big(|P_tf-\mu(f)|^2+ |\nn P_tf|^2\big)\le c_2 \e^{-\ll t} \mu\big(|f-\mu(f)|^2+|\nn f|^2\big),\ \ t\ge 0, f\in C_b^1(\R^d) $$ holds for some constants $c_2,\ll>0$, where and in the following,
$\mu(f):=\int f\d\mu$ for $f\in L^1(\mu)$. If  the gradient estimate $|\nn P_t f|^2 \le K(t)P_tf^2$ holds for some function $K: (0,\infty)\to (0,\infty)$, see  \cite{GW,WZ} for concrete estimates, we obtain the $L^2$-exponential convergence
\beq\label{L2} \Var_{\mu}(P_t f)\le c\e^{-\ll t} \Var_\mu(f),\ \ t\ge 0, f\in L^2(\mu)\end{equation} for some constants $c,\ll>0$, which has been derived in \cite{GS} using the idea of \cite{DMS}.
  See e.g. \cite{B,DMS,D,GM,GS,GW,W14,WZ} and references within for further results on exponential convergence and regularity estimates of $P_t$.

Recently,   Hu and Wang \cite{HW} prove the sub-exponential convergence by using the weak Poincar\'e inequality
\beq\label{WP} \Var_{\mu_1}(f)  \le \aa(r) \mu_1(|\nn f|^2) +r\|f\|_{\rm osc}^2,\ \ f\in C_b^1(\R^d)\end{equation}
for some decreasing function $\aa: (0,\infty)\to (0,\infty),$ where $\|f\|_{\rm osc}:={\rm ess}_\mu\sup f- {\rm ess}_\mu\inf f.$  According to \cite[Theorem 3.6]{HW}, \eqref{WP} implies
\beq\label{HW} \mu\big(|P_tf-\mu(f)|^2+ |\nn P_tf|^2\big)\le c_1\,\xi(t)   \big(\|f\|_\infty^2+\mu(|\nn f|^2)\big),\ \ t\ge 0, f\in C_b^1(\R^d) \end{equation}   for some constant $c_2>0$ and
$$\xi(t):= \inf\big\{s >0: \ t\ge -\aa(s) \log s\big\},\ \ t\ge 0.$$ Again, if the gradient estimate $|\nn P_t f|^2\le K(t)P_t f^2$ holds then this implies
\beq\label{WC} \Var_\mu(P_t f)\le  c_1\,\xi(t) \|f\|_{\rm osc}^2,\ \ t\ge 0, f\in L^\infty(\mu)\end{equation} for some constant $c_1>0$.
In particular, if $\aa$ is bounded so that \eqref{WP} reduces to \eqref{P} with $c_1= \|\aa\|_\infty$, we obtain the  exponential convergence as in the previous case.

In this paper we aim to introduce weak Poincar\'e inequalities to estimate the convergence rate for more general degenerate diffusion semigroups where $\mu_2$ is not necessarily a Gaussian measure.
Consider the following degenerate SDE for $(X_t,Y_t)$ on $\R^{d_1+d_2}=\R^{d_1}\times\R^{d_2}$, where $d_1,d_2\ge 1$ may be different:
\beq\label{E2} \beg{cases} \d X_t= Q(\nn^{(2)} V_2)(Y_t)\d t,\\
\d Y_t= \ss 2 \,\d B_t -\big(Q^*(\nn^{(1)} V_1) (X_t)+(\nn^{(2)} V_2)(Y_t)\big)\d t,\end{cases}\end{equation}
where $Q$ is a $d_1\times d_2$-matrix, $V_i\in C^2(\R^{d_i})$ such that $Z(V_i)<\infty, i=1,2,$ and $\nn^{(1)},\nn^{(2)}$ are the  gradient operators in components $x\in \R^{d_1}$ and $y\in \R^{d_2}$ respectively. It is easy to see that the generator of solutions to \eqref{E2} is dissipative in $L^2(\mu)$, where
$\mu:=\mu_1\times\mu_2$ for probability measures $\mu_i(\d x):= Z(V_i)^{-1} \e^{-V_i(x)}\d x$ on $\R^{d_i}, i=1,2;$  see the beginning of Section 3 for details. 

Since the coefficients of the SDE \eqref{E2} are locally Lipschitz continuous, for any initial point $z=(x,y)\in \R^{d_1+d_2}$, the SDE has a unique solution $(X_t^z, Y_t^z)$ up to life time $\zeta^z$.
Let $P_t$ be the associated (sub-) Markov semigroup, i.e.
$$P_t f(z)= \E\big[f(X_t^z, Y_t^z) 1_{\{t<\zeta^z\}}\big],\ \ f\in \B_b(\R^{d_1+d_2}), z\in \R^{d_1+d_2}, t\ge 0.$$

To ensure the non-explosion of the solution and the convergence of  the $L^2$-Markov  semigroup $P_t$   to $\mu$, we make the following assumption.

\beg{enumerate}\item[$(H)$] $QQ^*$ is invertible, there exists a constant $M>0$ such that
\beq\label{VV2} |(\nn^{(i)})^2 V_i|\le M (1+|\nn^{(i)} V_i|^{\tau_i}),\ \ i=1,2,
\end{equation}
for $\tau_1=1$ and some $1 \le \tau_2 <2$. Moreover,
$\mu_2(|\nn^{(2)} V_2|^4)<\infty$ and $V_2(y)=\Phi(|\sigma y-b|^2)$ for some invertible $d_2\times d_2$-matrix $\si$, $b\in \R^{d_2}$ and increasing function $\Phi\in C^3([0,\infty))$ such that
\beq\label{VV3} \sup_{r\ge 0} \Big|\Phi' (r)+2 r \Phi''(r)- \ff{2r \Phi'''(r) +(d_2+2)\Phi''(r)}{\Phi'(r)}\Big|  <\infty.\end{equation}
 \end{enumerate}

According to \cite[Theorem 3.1]{RW01},   there exist two decreasing functions  $\aa_1,\aa_2: (0,\infty)\to [1,\infty)$  such that the weak Poincar\'e inequality
\beq\label{WPi} \Var_{\mu_i}(f)  \le \aa_i(r) \mu_i(|\nn^{(i)} f|^2) +r\|f\|_{\rm osc}^2,\ \ f\in C_b^1(\R^{d_i}), r > 0,\end{equation}
holds for $i=1,2$.  We have the following result on the convergence rate of $P_t$ to $\mu$.

\beg{thm}\label{T1.1}  Let   $V_1$ and $V_2$ satisfy   $(H)$. Then the solution to $\eqref{E2}$ is non-explosive and $\mu$ is an invariant probability measure of  the associated  Markov semigroup $P_t$.    Moreover, there exist constants $c_1,c_2>0$ such that $\eqref{WC}$ holds for
\beq\label{WP*} \xi(t):= c_1 \inf\Big\{r >0: \ c_2 t \ge  \aa_1(r)^2\aa_2\Big(\ff r {\aa_1(r)^2}\Big)\log \ff 1 r\Big\},\end{equation}  which goes to $0$ as $t\to\infty.$
 \end{thm}

\paragraph{Remark 1.1.} (1) When $V_2(y)= \ff 1 2 |y|^2$ the measure $\mu_2$ reduces to the standard Gaussian measure as in \cite{HW}. In this case, we may repeat the argument in the proof of \cite[Theorem 3.6]{HW} to prove \eqref{HW} for
\beq\label{HP0} \xi(t)= \inf\Big\{r>0:\ c_2t \ge \aa_1(r)\log \ff 1 r\Big\},\ \ t>0,\end{equation}
and thus extend the main result in \cite{HW} to the case that $d_1\ne d_2$. Since in this case we have $\aa_2\equiv 1$, the convergence rate in Theorem \ref{T1.1} becomes
$$\xi(t)= \inf\Big\{r>0:\ c_2t \ge \aa_1(r)^2\log \ff 1 r\Big\},\ \ t>0,$$ which is in general worse than that in \eqref{HP0}. However, the argument in \cite{HW} heavily depends on the specific
$V_2(y)=\ff 1 2 |y|^2$ (or by linear  change of variables  $V_2(y)= |\si y-b|^2$ for some invertible $d_2\times d_2$-matrix $\si$ and $b\in \R^{d_2}$), and   is hard to extend to a general setting as in $(H)$. Nevertheless, we would hope to improve the convergence rate in Theorem \ref{T1.1} such that  \eqref{HP0} is covered for bounded $\aa_2$.

(2) Theorem \ref{T1.1} also applies to the following SDE for $(X_t,\bar Y_t)$ on $\R^{d_1+d_2}$ for some  invertible $d_2\times d_2$-matrix $\si$ and invertible $d_1\times d_1$-matrix $\bar Q\bar Q^*$:
\beq\label{E2'} \beg{cases} \d X_t= \bar Q(\nn^{(2)} V_2)(\bar Y_t)\d t,\\
\d \bar Y_t= \ss 2\, \si\d B_t -\big(\bar Q^*(\nn^{(1)} V_1) (X_t)+\si\si^*(\nn^{(2)} V_2)(\bar Y_t)\big)\d t.\end{cases}\end{equation}
Indeed, let $(X_t,Y_t)$ solve \eqref{E2} and let  $\bar Y_t= \si Y_t, \bar V_2(y)= V_2(\si^{-1}y)$. We have
$$(\nn^{(2)}\bar V_2)(y)= (\si^{-1})^* (\nn^{(2)} V_2)(\si^{-1} y),\ \ y\in \R^{d_2},$$ so that 
$$\d X_t = Q(\nn^{(2)}V_2)(Y_t)\d t= Q\si^* (\nn^{(2)}\bar V_2)(\bar Y_t),$$ and 
$$\d \bar Y_t= \ss 2\, \si \d B_t - \big(\si Q^* (\nn^{(1)}V_1)(X_t) + \si\si^* (\nn^{(2)} \bar V_2)(\bar Y_t)\big)\d t.$$
Letting $\bar Q= Q\si^*$, we see that the SDE \eqref{E2} is equivalent  to  \eqref{E2'}.  

 \

 To illustrate Theorem \ref{T1.1}, we consider the following   example with some concrete convergence rates of $P_t$.

 \paragraph{Example 1.1.}  We write $f\sim g$ for  real functions $f$ and $g$ on $\R^d$ if $f-g\in C_b^2(\R^d).$

 \ \newline
 {\bf (A)} Let $V_1(x)\sim k(1+|x|^2)^{\ff \dd 2 }$  for some constants $k,\dd>0$.
 \beg{enumerate}
\item[$(A_1)$] When $V_2(y)=   \kk (1 + |y|^2 )^{\ff \vv 2}$ for some constants $ \kk,\vv>0$,  \eqref{WC} holds with
    $$\xi(t)= c_1  \exp\Big(- c_2 t^{\ff{\vv\dd}{\vv\dd + 8\vv(1-\dd)^+ + 4\dd(1-\vv)^+}}\Big),\ \ t\ge 0,$$ for some constants $c_1,c_2>0$. If, in particular,
    $\dd,\vv\ge 1$  then $P_t$ converges to $\mu$ exponentially fast.
\item[$(A_2)$] When  $V_2(y)=  \ff{d+p}2 \log (1+|y|^2)$  for some constant  $  p >0$,
  \eqref{WC} holds with
    $$\xi(t)= c (1+t)^{-\ff 1 {\theta(p)}}  \big(\log(\e+t)\big)^{\ff{8(\theta(p)+1)(1-\dd)^+ +\dd}{\theta(p)\dd}} $$    for some constant  $c >0$ and
   $$\theta(p):= \ff{d+p+2}{p}\land \ff{4p+4+2d}{(p^2-4-2d-2p)^+}.$$
\item[$(A_3)$] When $V_2(y)=  \ff{d}2 \log (1+|y|^2) + p \log\log(\e +|y|^2) $ for some constants    $p>1$,  \eqref{WC} holds with
    $$\xi(t)= c_1 \big(\log(\e+t)\big)^{1-p} \cdot\big(\log\log(\e^2+t) \big)^{\ff{8(1-\dd)^+}\dd}$$   for some constants $c>0$.
 \end{enumerate}
\ \newline
  {\bf (B)}  Let $V_1(x)\sim \ff{d+q}2 \log(1+|x|^2)$  for some  $q>0.$
 \beg{enumerate}
 \item[$(B_1)$]    When $V_2(y)=   k (1 + |y|^2 )^{\ff \vv 2}$ for some constants $ k,\vv >0$,  \eqref{WC} holds with
    $$\xi(t)= c  (1+t)^{-\ff 1 {2\theta(q)}} \big(\log(\e +t)\big)^{\ff {4(1-\vv)^++\vv}{2\vv\theta(q)}} $$ for some constant $c>0.$
 \item[$(B_2)$]  When $V_2(y)=  \ff{p+d}2\log(1+|y|^2)$ for some constant $ p >0$,  \eqref{WC} holds with
    $$\xi(t)= c (1+ t)^{-\ff 1 {2\theta(q)+\theta(p)+2\theta(p)\theta(q)}}\big(\log(\e+t)\big)^{\ff 1{2\theta(q)+\theta(p)+2\theta(p)\theta(q)}} $$ for some constant $c>0$.
 \item[$(B_3)$]  When $V_2(y)=     \ff{d}2 \log (1+|y|^2) + p \log\log(\e +|y|^2) $ for some constant  $p>1$,   \eqref{WC} holds with
    $$\xi(t)= c   \big(\log (\e +t)\big)^{-\ff{p-1}{1+2\theta(q)}}  $$ for some constant $c>0$.  \end{enumerate}
\ \newline
 {\bf (C)}  Let $V_1(x)\sim \ff{d}2 \log(1+|x|^2) +q \log\log (\e +|x|^2)$  for some  $q>0.$
 \beg{enumerate}
 \item[$(C_1)$]    When $V_2(y)=   k (1 + |y|^2 )^{\ff \vv 2}$ for some constant $ k >0$ and $\vv>0$,  or \newline
 $V_2(y)=  \ff{p+d}2\log(1+|y|^2)$ for some constant $ p >0$, \eqref{WC} holds with
    $$\xi(t)= c\big( \log (\e+t)\big)^{-(q-1)}$$ for some constant $c>0.$
 \item[$(C_2)$]   When $V_2(y)=     \ff{d}2 \log (1+|y|^2) + p \log\log(\e +|y|^2) $ for some constant  $p>1$,   \eqref{WC} holds with
    $$\xi(t)= c   \big(\log \log (\e^2 +t)\big)^{-(q-1)}  $$ for some constant $c>0$.  \end{enumerate}

 \

 In the next section we present a general  result on the weak hypocoercivity for $C_0$-semigroups on Hilbert spaces, which  is then used in Section 3 to prove Theorem \ref{T1.1} and Example 1.1.

 \section{A general framework}

 Let $(\H, \<\cdot, \cdot\>, \|\cdot\|)$ be a separable Hilbert space, let $(L,\D(L))$ be a densely defined linear operator generating a  $C_0$- contraction semigroup $P_t=\e^{tL}$. We aim to investigate the decay rate of $P_t$ of type
 \beq\label{2*1} \|P_tf\|^2\le \xi(t) \big(\|f\|^2+\Psi(f)\big),\ \ t\ge 0, f\in \D(L),\end{equation} where $\xi$ is a decreasing function with $\xi(\infty):=\lim_{t\to\infty}\xi(t)=0$, and $\Psi: \H\to [0,\infty]$ is a functional such that the set $\{f\in \H: \Psi(f)<\infty\}$ is dense in $\H$. 

 \subsection{Main result}

Following the line of e.g. \cite{DMS,GS}, we assume that $L$ decomposes into  symmetric and antisymmetric part:
$$L= S-A\ \ \text{on} \ \D,$$ where $\D$ is a core of $(L,\D(L))$,   $S$ is symmetric and    $A$ is    antisymmetric. Then both $(S,\D)$ and $(A,\D)$ are closable in $\H$. Let $(S,\D(S))$ and $(A,\D(A))$ be their closures.
These two operators are   linked to  the orthogonal decomposition $\H=\H_1\oplus \H_2$ in the following assumptions, where
$$\pi_i: \H\to \H_i,\ \ i=1,2,$$ are the orthogonal projections.

\beg{enumerate} \item[$(H1)$] $\H_1\subset {\mathcal N} (S):=\{f\in\D(S):\ Sf=0\}$; that is, $\H_1\subset \D(S)$ (hence, $\pi_2 \D\subset \D(S)$ due to $\D\subset \D(S)$) and $S\pi_1=0.$
\item[$(H2)$] $\pi_1\D\subset \D(A)$ (hence, also $\pi_2 \D\subset \D(A)$ due to $\D\subset \D(A)$)  and $\pi_1 A\pi_1|_\D=0.$
\end{enumerate}

Since $(A,\D(A))$ is closed, antisymmetric and $\pi_1\D \subset \D(A)$, $(\pi_1A,\D(A))$ is closable. Denote the closure by $(\pi_1A, \D(\pi_1A))$. By $(H2)$, $A\pi_1$ is well defined on $\D$, and by the antisymmetry of $A$,
$$(A\pi_1)^*= \pi_1 A^*= -\pi_1 A\ \text{holds\ on\ } \D.$$
Then $A\pi_1$ with domain $\D(A\pi_1):=\{f\in \H: \pi_1f\in \D(A)\}$ is a densely defined closed operator with adjoint $(-\pi_1 A, \D(\pi_1 A))$.
By von Neumann's theorem, see e.g.~\cite[Theorem 5.1.9]{Ped},   the operators $G := (A\pi_1)^*A\pi_1$ and $I+(A\pi_1)^*A\pi_1$ with domain
\begin{align*}
\D(G) := \D((A\pi_1)^*(A\pi_1))=\big\{f \in \D(A\pi_1): ~A\pi_1 f \in \D((A\pi_1)^*)\,\big\}
\end{align*}
are are self-adjoint. Furthermore, the latter one is injective and surjective (with range equal to $\H$) and admits a bounded linear inverse. We define the operator $B$ with domain $\D(B)=\D((A\pi_1)^*)$ via
\begin{align} \label{1.1}
B:=(I+(A\pi_1)^*A\pi_1)^{-1}(A\pi_1)^*.
\end{align}
Then   $B^*=A\pi_1(I+(A\pi_1)^*A\pi_1)^{-1}$ defined on $\D(B^*)=\H$ is closed and   bounded. Consequently,   $(B,D((A\pi_1)^*))$ is also bounded and has a unique extension to a bounded linear operator $(B,\H)$. By e.g.  \cite[Theorem 5.1.9]{Ped}, we have
\begin{align*}
B=(A\pi_1)^*(I+A\pi_1(A\pi_1)^*)^{-1}.
\end{align*} Consequently, $\|B\|\le 1$ and $\pi_1B=B$.
%Denote by $(G, \D(G))$ the self-adjoint operator $((A\pi_1)^*(A\pi_1), \D((A\pi_1)^*(A\pi_1)))$.

We shall need the following two more assumptions.

\beg{enumerate}
\item[$(H3)$]We assume $\D \subset \D(G)$.
%and $\D$ is a core of $(G, \D(G))$.
Furthermore, there exists a constant $N\ge 1$ such that
$$\big\<BS\pi_2f, \pi_1 f\big\>\le \frac{N}{2}\|\pi_1f\|\cdot\|\pi_2f\|,\ \ -\big\<BA\pi_2f, \pi_1 f\big\>\le \frac{N}{2}\|\pi_1f\|\cdot\|\pi_2f\|,\ \ f\in \D.$$

\item[$(H4)$] For any $f\in \D(L)$ there exists a sequence $\{f_n\}_{n\ge 1}\subset \D$ such that $f_n\to f$ in $\H$ and
$$\limsup_{n\to \infty} \<-Lf_n, f_n\>\le \<-L f, f\>, \ \ \limsup_{n\to\infty} \Psi(f_n)\le \Psi(f).$$ %Moreover, for any $f \in \D(A\pi_1)$
%there exists a sequence $\{f_n\}_{n\ge 1}\subset \D$ such that $f_n\to f$ in $\H$ and
%$$\limsup_{n\to \infty} \<A\pi_1 f_n, A\pi_1 f_n\>\le \<\A\pi_1 f, A\pi_1 f\>, \ \ \limsup_{n\to\infty} \Psi(f_n)\le \Psi(f).$$
\end{enumerate}

\beg{thm} \label{T2.1} Assume $(H1)$-$(H4)$ and let $\Psi$ satisfy
\beq\label{phi} \Psi(P_tf)\le \Psi(f),\ \ \Psi(\e^{-tG}f)\le \Psi(f),\ \ \Psi(\pi_1f)\le \Psi(f),\ \ f\in \H.\end{equation} If the weak Poincar\'e inequalities
\beq\label{WP1} \|\pi_1 f\|^2\le \aa_1(r) \|A\pi_1f\|^2 + r\Psi(\pi_1 f),\ \ r>0, f\in \D(A\pi_1), \end{equation} and
\beq\label{WP2} \|\pi_2 f\|^2\le \aa_2(r) \<-S f,  f\>  + r\Psi(f),\ \ r>0, f\in \D \end{equation} hold for some  decreasing
functions $\aa_i: (0,\infty)\to [1,\infty), i=1,2,$  then there exist constants $c_1,c_2>0$ such that  $\eqref{2*1}$ holds for
\beq\label{WPG*} \xi(t):= c_1 \inf\Big\{r >0: \ c_2 t \ge  \aa_1(r)^2\aa_2\Big(\ff r{\aa_1(r)^2}\Big) \log \ff 1 r \Big\},\end{equation}  which goes to $0$ as $t\to\infty.$
\end{thm}

\subsection{Preparations}

 \beg{lem}\label{L2.1} Under  $(H1)$-$(H3)$, we have
 \beq\label{BA1} \|Bf\| \le \ff 1 2 \|\pi_2 f\|,\ \ f\in \H,\end{equation}
 \beq\label{BA2} \|ABf\| \le   \|\pi_2 f\|,\ \ f\in \D,\end{equation}
 \beq\label{BA3} |\<Bf, Lf\>| \le \|\pi_2f\|\cdot \|f\|,\ \ f\in \D(L),\end{equation}
 \beq\label{BA4} \<BLf, f\>\le N\|\pi_1f\|\cdot\|\pi_2 f\|-\big\<(1+G)^{-1} G\pi_1f, \pi_1f\big\>,\ \ f\in\D(L).\end{equation}\end{lem}

\beg{proof} Let $f\in \D$ and $g= Bf.$ By \eqref{1.1}, $\pi_1A^*\pi_1f=-\pi_1A\pi_1f=0$ and $\pi_2 f\in \D(A)$ (see $(H2)$), we have
\beq\label{ABD}\beg{split} &\|g\|^2+\|A\pi_1g\|^2= \<g+(A\pi_1)^*A\pi_1 g, g\> = \<(A\pi_1)^*f, g\>\\
&= \<(A\pi_1)^*\pi_2f, g\>=\<\pi_2 f, A\pi_1 g\> \le \|\pi_2 f\|\cdot \| A\pi_1 g\|.\end{split}\end{equation}
Combining this with $$\|\pi_2 f\|\cdot \| A\pi_1 g\|\le \ff 1 4 \|\pi_2 f\|^2 +\|A\pi_1 g\|^2,$$ we obtain \eqref{BA1} for $f\in\D$, and hence for all $f\in\H$ since $\D$ is dense in $\H$ and the operators $B, \pi_2$ are bounded.

Next, combining \eqref{ABD} with $\pi_1 B=B$ and  $$\|\pi_2 f\|\cdot \| A\pi_1 g\|\le \ff 1 2 \|\pi_2 f\|^2 +\ff 1 2 \|A\pi_1 g\|^2,$$ we obtain
$$\|ABf\|^2=\|A\pi_1Bf\|^2=  \|A\pi_1 g\|^2 \le \|\pi_2 f\|^2,\ \ f\in\D,$$ which is equivalent to \eqref{BA2}.

 Moreover, by the symmetry of $S,$ antisymmetry of $A$,
$S\pi_1=0,$ and $ B= \pi_1 B$,  we obtain from \eqref{BA2} that for any $f\in \D$,
$$|\<Bf, Lf\>|= |\<Bf, -A f\> |= |\<ABf, f\>| \le \|\pi_2 f\|\cdot \|f\|.$$ Since $\D$ is dense in $\D(L)$ and $B$ is bounded, this implies \eqref{BA3}.

Finally, by $\pi_1B=B$, $S\pi_1=0$, the definition of $B$ and $(H3)$, for $f\in\D$ we have
\beg{align*} &\<BLf, f\> = \<BLf, \pi_1 f\>= \<BS f, \pi_1 f\>-\<BA f, \pi_1 f\> \\
&=\<BS \pi_2 f, \pi_1 f\>\\-\<BA \pi_1 f, \pi_1 f\> - \<BA \pi_2 f, \pi_1 f\>\\
&\le N\|\pi_1f\|\cdot\|\pi_2f\|-\big\<(1+G)^{-1} G \pi_1f, \pi_1f\big\>.
\end{align*}
By the boundedness of $(1+G)^{-1} G$ and that $\D$ is dense in $\D(L)$, this implies
\eqref{BA4}.
\end{proof}

Next, we need the following  result on weak Poincar\'e inequality for subordinated operators. Let $\nu$ be a L\'evy measure on $[0,\infty)$ such that $\int_0^\infty (r\land 1)\nu(\d r)<\infty$, then
$$\phi_\nu(s):= \int_0^\infty\big(1-\e^{-sr}\big)\nu(\d r),\ \ s\ge 0$$ is a Bernstein function. Let $(S_0,\D(S_0))$ be a non-negative definite self-adjoint operator. We intend to establish the weak Poincar\'e inequality for the form $\<\phi_\nu(S_0)f,f\>$ in terms of that for $\<S_0f,f\>$. The Nash type and  super Poincar\'e  inequalities have already been investigated in \cite{BM,SW}. Recently, sub-exponential decay for   subordinated semigroups was studied in \cite{DS}, where $\phi_\nu$ is assumed to satisfy
$$\liminf_{s\to\infty} \ff{\phi_\nu(s)}{\log s}>0.$$ However, this condition  excludes  $\phi_\nu(s):= \ff s {1+s}$
which is indeed what we need   in the proof of Theorem \ref{T2.1}.

\beg{lem}\label{L2.2} Let $(A_0,\D(A_0))$ be a densely defined closed linear operator on a separable Hilbert space $\H_0$. Let $P_t^0$ be the $C_0$-contraction semigroup   generated by the self-adjoint operator  $-A_0^*A_0$ with domain $\D(A_0^*A_0):= \{f\in \D(A_0): A_0f\in \D(A_0^*)\}$.   If the weak Poincar\'e inequality
\beq\label{WP0} \|f\|^2\le \aa(r) \|A_0 f\|^2 +r\Psi_0(f),\ \ r>0, f\in \D(A_0)\end{equation} holds for some decreasing $\aa: (0,\infty)\to (0,\infty)$,  where $\Psi_0: \H_0\to [0,\infty]$ satisfies
\beq\label{PHI} \Psi_0(P_t^0f)\le \Psi_0(f),\ \ t\ge 0, f\in \D(A_0),\end{equation} then
$$\|f\|^2\le \bigg(\int_0^\infty \big(1-\e^{-\ff s{\aa(r)}}\big)\nu(\d s) \bigg)^{-1}  \big\|(\phi_\nu(A_0^*A_0))^{1/2} f \big\|^2 + r \Psi(f),\ \ r>0, f\in \D(A_0).$$ In particular, for $\nu(\d s)= \e^{-s}\d s$ such that $\phi_\nu(s)=\ff s {1+s}$, we have
$$\|f\|^2\le  \big(1+\aa(r)\big) \big\<(1+A_0^*A_0)^{-1}A_0^*A_0 f,f\big\> + r \Psi(f),\ \ r>0, f\in\D(A_0).$$
\end{lem}

\beg{proof}  Since $\D((A_0^*A_0)^{1/2})=\D(A_0)$, we have $\D(\{\phi_\nu(A_0^*A_0)\}^{1/2})\supset \D(A_0)$.  By \eqref{WP0} and \eqref{PHI}, for any $f\in \D(A_0)$,
 $$\ff{\d }{\d t} \|P_t^0f\|^2=-2\|A_0  P_t^0f\|^2 \le -\ff 2 {\aa(r)}\|P_t^0f\|^2 +\ff{2r}{\aa(r)}\Psi(f),\ \ t\ge 0, r>0,$$ because
 $P_t^0$ leaves $\D(A_0)$ invariant. Then Gronwall's lemma gives  \beq\label{ET2}\|P_t^0f\|^2\le \e^{-\ff{2t}{\aa(r)}} \|f\|^2 + r\Psi(f)(1-\e^{-\ff{2t}{\aa(r)}}),\ \ r>0, t\ge 0.\end{equation} Therefore,
\beg{align*} &\Big\|\big(\phi_\nu(A_0^*A_0)\big)^{1/2}  f \Big\|^2 = \int_0^\infty\<f-P_s^0f, f\>\nu(\d s)= \int_0^\infty\big(\|f\|^2-\|P_{s/2}^0f\|^2\big)\nu(\d s)\\
&\ge  \int_0^\infty \Big(\|f\|^2-\e^{-\ff{s}{\aa(r)}} \|f\|^2 - r\Psi(f)(1-\e^{-\ff{s}{\aa(r)}})\Big)\nu(\d s)\\
&=  \big(\|f\|^2-r\Psi(r)\big)\int_0^\infty \big(1-\e^{-\ff{s}{\aa(r)}}\big)\nu(\d s),\ \ r>0.\end{align*}
This implies  the desired inequality.
\end{proof}

In the proof of Theorem \ref{T1.1} (see Section 3 below), to verify $(H3)$ we check the following two inequalities:
\beq\label{HH3}\beg{split}  & \big\<BS\pi_2f, \pi_1 f\big\>\le N\|\pi_1f\|\cdot\|\pi_2f\|,\\
&\big\<BA\pi_2f, \pi_1 f\big\>\le N\|\pi_1f\|\cdot\|\pi_2f\|,\ \ f\in\D.\end{split} \end{equation} The first inequality is easy to check there, see Section 3, the first part in the proof of $(H3)$.  To verify the second, we present below a sufficient condition provided in \cite[Prop.~2.15]{GS}.

\begin{prp} \label{P1.14}
Assume that $(-G, \D)$ is essentially m-dissipative (equivalently, essentially self-adjoint). If    there exists constant $N\in (0,\infty)$ such that
\beq \label{BA*}
\|(BA)^*g\| \leq N\,\|g\| ~~\mbox{ for all }~ g=(I+G)f,~f \in \D,
\end{equation}
then $$\big|\big\<BA\pi_2f, \pi_1 f\big\>\big|\le N\|\pi_1f\|\cdot\|\pi_2f\|,\ \ f\in\D.$$
\end{prp}

\subsection{Proof of Theorem \ref{T2.1}}

\beg{proof}
For any $\vv\in [0,1)$, let
$$I_\vv(f)= \ff 1 2 \|f\|^2+\vv\<Bf, f\>,\ \ f\in\H.$$ By \eqref{BA1}, we have
\beq\label{BA5} \ff{1-\vv} 2 \|f\|^2\le I_\vv(f)\le \ff{1+\vv} 2 \|f\|^2,\ \ f\in\H.\end{equation}
Now, let $f\in \D$ and $f_t= P_t f$ for $t\ge 0$. We have
\beq\label{3.1} \ff{\d}{\d t} I_\vv(f_t)= \<Lf_t, f_t\> +\vv \<BLf_t, f_t\>+ \vv \<Bf_t, Lf_t\>.\end{equation}
By   \eqref{WP2} and $\<-Lg,g\>=\<-Sg,g\>$ for $g\in \D$, we obtain
$$ \<Lg, g\>  \le -\ff{\|\pi_2g\|^2}{\aa_2(r_2)} + \ff{r_2\Psi(g)}{\aa_2(r_2)},\ \  g\in\D,r_2>0.$$
Since $f_t\in \D(L)$, combining this with  $(H4)$ and  \eqref{phi}, we arrive at
\beq\label{3.2} \<Lf_t, f_t\>  \le -\ff{\|\pi_2 f_t\|^2}{\aa_2(r_2)} + \ff{r_2\Psi(f_t)}{\aa_2(r_2)}\le -\ff{\|\pi_2 f_t\|^2}{\aa_2(r_2)} + \ff{r_2\Psi(f)}{\aa_2(r_2)},\ \ t,r_2>0.\end{equation}
Next,  applying Lemma \ref{L2.2} with $\H_0=\H_1, A_0= ((A\pi_1)^*(A\pi_1))^{1/2}|_{\H_1}$ and $\Psi_0=\Psi|_{\H_1}$ such that condition \eqref{PHI} follows from \eqref{phi}, we see that \eqref{WP1} implies
$$-\big\<(I+(A\pi_1)^*(A\pi_1))^{-1} (A\pi_1)^*A\pi_1 f, \pi_1 f\big\> \le   -\ff{\|\pi_1 f\|^2}{\aa_1(r_1)+1} +\ff{r_1\Psi(\pi_1 f)}{\aa_1(r_1)+1},\ \ r>0, f \in \D(A\pi_1).$$ Since the operator $(I+(A\pi_1)^*(A\pi_1))^{-1} (A\pi_1)^*A\pi_1 $ is bounded, $\D(A\pi_1)\supset \D$ due to $(H2)$,   and by $(H4)$ for any $g\in \D(L)$ we may find a sequence $g_n\in \D$ such that $g_n\to g$ in $\H$ and $\limsup_{n\to\infty} \Psi(g_n)\le \Psi(g)$, this inequality holds for all $g\in \D(L).$
 Combining this with \eqref{BA4} and \eqref{phi},  we obtain
\beq\label{3.3}\beg{split} \<BLf_t, f_t\>&\le N\|\pi_1 f_t\|\cdot\|\pi_2 f_t\| - \big\<(I+(A\pi_1)^*(A\pi_1))^{-1} (A\pi_1)^*A\pi_1 f_t, \pi_1 f_t\big\>\\
&\le N\|\pi_1 f_t\|\cdot\|\pi_2 f_t\| -\ff{\|\pi_1 f_t\|^2}{\aa_1(r_1)+1} +\ff{r_1\Psi(f)}{\aa_1(r_1)+1},\ \ t,r_1>0.\end{split}\end{equation}
Substituting \eqref{BA3}, \eqref{3.2} and \eqref{3.3} into \eqref{3.1}, we arrive at
\beg{align*} \ff{\d}{\d t} I_\vv(f_t)\le &-\Big(\ff{\|\pi_2 f_t\|^2}{\aa_2(r_2)} +\ff{\vv \|\pi_1f_t\|^2}{\aa_1(r_1)+1}\Big)   +\vv\big(N\|\pi_1 f_t\|\cdot \|\pi_2 f_t\|+ \|\pi_2 f_t\|\cdot\|f_t\|\big)\\
&+\Psi(f)\Big(\ff{r_2}{\aa_2(r_2)}+ \ff{\vv r_1}{\aa_1(r_1)+1}\Big),\ \ t\ge 0, f\in\D.\end{align*}
Combining this with
\beg{align*} &\vv N \|\pi_1 f_t\|\cdot \|\pi_2 f_t\|\le \ff{\vv \|\pi_1 f_t\|^2}{2(\aa_1(r_1)+1)} +\ff{\vv N^2 (\aa_1(r_1)+1)\|\pi_2 f_t\|^2}{2},\\
&\vv\|\pi_2 f_t\|\cdot \|f_t\|\le \ff{\|\pi_2 f_t\|^2}{2\aa_2(r_2)} + \ff{\vv^2 \aa_2(r_2)\|f_t\|^2} 2,\end{align*}
we obtain
\beq\label{II} \beg{split} \ff{\d}{\d t} I_\vv(f_t)&\le -\Big(\ff 1 {2\aa_2(r_2)} - \ff{\vv N^2(\aa_1(r_1)+1)}2\Big)  \|\pi_2 f_t\|^2 -\ff{\vv  \|\pi_1f_t\|^2}{2(\aa_1(r_1)+1)} \\
&\quad  +\ff{\vv^2\aa_2(r_2)\|f_t\|^2} 2 + \Psi(f)\Big(\ff{r_2}{\aa_2(r_2)}+ \ff{\vv r_1}{\aa_1(r_1)+1}\Big),\ \ t\ge 0, f\in\D.\end{split}\end{equation}
Taking \beq\label{*V} \vv= \ff 1 {2N^2(\aa_1(r_1)+1)\aa_2(r_2)}\le \ff 1 2\end{equation} since $N,\aa_2\ge 1$, we have
\beg{align*} &\ff 1 {2\aa_2(r_2)}- \ff{\vv N^2(\aa_1(r_1)+1)}{2} \ge \ff 1 {4\aa_2(r_2)},\\
&\ff{1}{4\aa_2(r_2)}\land \ff{\vv}{2(\aa_1(r_1)+1)}\ge \vv^2 \aa_2(r_2).\end{align*}
Then \eqref{II} implies
$$\ff{\d}{\d t} I_\vv(f_t) \le -\ff{ \|f_t\|^2} {8N^4\aa_2(r_2)(\aa_1(r_1)+1)^2}   +   \Psi(f)\Big(\ff{r_2}{\aa_2(r_2)}+ \ff{r_1}{2N^2 \aa_2(r_2)(\aa_1(r_1)+1)^2}\Big).$$
Since $\vv\le \ff 1 2$, by \eqref{BA5} we have $\|f_t\|^2\ge \ff 4 3 I_\vv(f_t)$, so that
$$\ff{\d}{\d t} I_\vv(f_t) \le -\ff{ I_\vv(f_t)  } {6N^4\aa_2(r_2)(\aa_1(r_1)+1)^2}  +  \Psi(f)\Big(\ff{r_2}{\aa_2(r_2)}+ \ff{r_1}{2N^2 \aa_2(r_2)(\aa_1(r_1)+1)^2}\Big).$$
By Gronwall's lemma and \eqref{*V}, we arrive at
\beg{align*} I_\vv(f_t)&\le \exp\Big[-\ff {t}{6N^4\aa_2(r_2)(\aa_1(r_1)+1)^2}\Big] I_\vv(f)+  \Psi(f)\Big(3N^2r_1+6N^4r_2(\aa_1(r_1)+1)^2\Big).  \end{align*}
Taking $r_1=r, r_2=\ff r {\aa_1(r)^2}$, using   \eqref{BA5} for $\vv\in (0,\ff 1 2)$ and  that $\aa_1(r)\ge 1$, obtain
$$\|f_t\|^2\le c_1\exp\Big[-\ff{c_2t}{\aa_1(r)^2\aa_2(\ff r{\aa_1(r)^2})}\Big] \|f\|^2 + c_1 r  \Psi(f),\ \ r>0, f\in\D, t\ge 0.$$  Consequently,  for any $r>0$  such that $c_2t \ge  \aa_1(r)^2\aa_2 (\ff r{\aa_1(r)^2})\log \ff 1 r,$ we have
$$\|f_t\|^2\le c_1r \big(\|f\|^2+ \Psi(f)\big).$$ Therefore,  \eqref{2*1} with $\xi(t)$ in \eqref{WPG*} holds for $f\in\D$. By $(H4)$,  it holds  for all $f\in\D(L).$ Then the proof is finished. \end{proof}

\section{Proof of Theorem \ref{T1.1}}

We first embed $P_t$ in the framework of Section 2. By shifting the second variable $y$, in $(H)$ we may and do  take $b=0$, i.e. $V_2(y)=\Phi(|\si y|^2),$ for some invertible $d_2\times d_2$-matrix $\si$. Since we may move $\sigma$ from the potential $V_2$ to the symmetric part of the generator $L$ corresponding to the solution of \eqref{E2} and the matrix $Q$ as described in Remark 1.1(2), we only have to consider the case $V_2(y)=\Phi(|y|^2)$. Thus
\beq\label{QP} \nn^{(2)}V_2(y)= 2\Phi'(|y|^2)\,y.\end{equation}
Let
\beq\label{Mu} \mu=\mu_1\times\mu_2,\ \ \text{where}\  \mu_i(\d x_i):= Z(V_i)^{-1}\e^{-V_i(x_i)}\d x_i\ \text{on}\ \R^{d_1},\  i=1,2.\end{equation}
By It\^o's formula, the generator $L$ for the solution to \eqref{E2} has the decomposition
$$L= S-A,$$ where \beg{align*} S:=&\DD^{(2)} - \<(\nn^{(2)} V_2), \nn^{(2)}\cdot\big\>=\sum_{i=1}^{d_2} \big(\pp_{y_i}^2 - (\pp_{y_i}V_2)\pp_{y_i}\big),\\
 A:= &\big\<Q^*(\nn^{(1)} V_1), \nn^{(2)}\cdot\big\> - \big\<Q(\nn^{(2)} V_2), \nn^{(1)}\cdot\big\>\\
 = &\sum_{i=1}^{d_1}\sum_{j=1}^{d_2}  Q_{ij} \big((\pp_{x_i}V_1)\pp_{y_j}  - (\pp_{y_j} V_2)\pp_{x_i}\big).\end{align*}
Since above we moved $\sigma$ from the potential $V_2$ to the symmetric part of $L$ and to the matrix $Q$, instead of $S$ and $Q$ we should consider  
\beg{align*} \sum_{i,j=1}^{d_2} (\sigma\sigma^*)_{ij}\big(\pp_{y_i}\pp_{y_j} - (\pp_{y_i}V_2)\pp_{y_j}\big) \quad \mbox{and} \quad  Q\sigma^*, \end{align*}
respectively. But, because $\sigma\sigma^*$ is a constant, symmetric, invertible matrix, without loss of generality we may take $\sigma$ equal to the identity matrix. The considerations below easily generalize to general $\sigma$, but are easier to follow for $\sigma$ being the identity matrix.

Let $\nn= (\nn^{(1)}, \nn^{(2)})$ be the gradient operator on $\R^{d_1+d_2}$, and denote
$$C_c^\infty(\R^{d_1+d_2})=\big\{f\in C^\infty (\R^{d_1+d_2}): \nn f\ \text{has\ compact\ support}\big\}.$$   The integration by parts formula implies that
$(S, C_c^\infty(\R^{d_1+d_2}))$ is symmetric and non-positive definite in $L^2(\mu)$ while   $(A,C_c^\infty(\R^{d_1+d_2}))$ is antisymmetric in $L^2(\mu)$. Consequently,
   $L^*:=L+2A= S+A$ satisfies  $$\mu(fLg)=\mu(gL^*f),\ \ f,g\in C_c^\infty(\R^{d_1+d_2}).$$
  Therefore,  $(L, C_c^\infty(\R^{d_1+d_2})$ is dissipative and, in particular, closable in $L^2(\mu)$. Let $(L,\D(L))$ denote the closure.
Then the first assertion of Theorem \ref{T1.1} is implied by the following proposition.

\beg{prp}\label{P3.1} Under assumption $(H)$,  the operator $(L, C_c^\infty(\R^{d_1+d_2}))$ is essentially m-dissipative in $L^2(\mu)$, and the $C_0$-contraction semigroup $T_t$ generated by the closure  coincides with $P_t$ in $L^2(\mu)$. Consequently, the solution to $\eqref{E2}$ is non-explosive and $\mu$ is an invariant probability measure of $P_t$. \end{prp}

\beg{proof} In \cite[Theorem 3.10]{nonne16} under even weaker assumptions as in $(H)$, essential m-dissipati\-vity of $(L, C_c^\infty(\R^{d_1+d_2}))$ in $L^2(\mu)$ is shown. In the proof condition \eqref{VV2} for $i=2$ is used. Hence the closure  $(L, \D(L))$ generates a $C_0$-contraction semigroup $T_t$. Then $\mu(Lf)=0$ for $f\in \D(L)$ implies that
$$\pp_t \mu(T_t f)= \mu(LT_t f)=0,\ \  t\ge 0, f\in \D(L),$$ so that    $\mu$ is an invariant probability measure of $T_t$.  On the other hand, according to \cite[Theorem 1.1 and Proposition 1.4]{BBR}
(see also \cite[Theorem 3.17 and Remark 3.18]{CG10}), for $\mu$-a.e.~starting point $z=(x,y) \in \R^{d_1+d_2}$ there is a law ${\mathbb P}^z$ on the space of $\R^{d_1+d_2}$-valued continuous functions such that $(X_t,Y_t)_{t\ge 0}$ is a weak solution to \eqref{E2} and for any distribution $\nu(\d z)= \rr(z)\mu(\d z)$ with a probability density $\rr$,
$$\mu(\rr T_t f)=\int_{\R^{d_1+d_2}}\E^z\big[f(X_t,Y_t)]\,\nu(dz),\ \ t\ge 0, f\in \B_b(\R^{d_1+d_2}).$$ By the uniqueness of  the SDE \eqref{E2}, we have for $\mu$-a.e.~$z \in \R^{d_1+d_2}$:
$$P_t f(x,y)= \E^z\big[f(X_t,Y_t)],\ \ t\ge 0, f\in \B_b(\R^{d_1+d_2}).$$  Therefore, $\mu(\rr P_tf)=\mu(\rr T_tf)$  holds for any $\rr\in L^1(\mu),  t\ge 0$ and $f\in \B_b(\R^{d_1+d_2})$,
and hence, $P_t$ is a $\mu$-version of $T_t$. Consequently, $\mu$ is an invariant probability measure of $P_t$. Since $P_t 1\le 1$, this implies that $P_t 1=1,\ \mu$-a.e. Since the coefficients of the SDE is at least $C^1$-smooth, the semigroup $P_t$ is Feller so that $P_t 1$ is continuous. Therefore, $P_t 1(z)=1$ holds for all $z\in \R^{d_1+d_2}$, i.e. the solution to \eqref{E2} is non-explosive. \end{proof}

Now,  to prove the second assertion in Theorem \ref{T1.1} using
 Theorem \ref{T2.1},    we take
$$\H=\{f\in L^2(\mu):\ \mu(f)=0\},\ \ \H_1=\{f\in \H:\ f(x,y) \text{\ does \ not\ depend\ on\ } y\}.$$
Then
\beq\label{5.0} (\pi_1f)(x,y)=\pi_1 f(x):= \int_{\R^{d_2}} f(x,y)\mu_2(\d y),\ \ f\in \H.\end{equation}
Let $$\D= \H\cap C_c^\infty(\R^{d_1+d_2})=\Big\{f\in  C_c^\infty(\R^{d_1+d_2}):\ \mu(f)=0\Big\}.$$  Let $(L,\D(L)), (S,\D(S))$ and $(A,\D(A))$ be the closures in $\H$ of
$(L, \D), (S, \D)$ and $(A,\D)$ respectively. Since the closure of $(L, C_c^\infty(\R^{d_1+d_2}))$ in $L^2(\mu)$ generates a strongly continuous contraction semigroup, see Proposition \ref{P3.1}, we have $L^2(\mu) = \overline{{\mathcal R}(L)} \oplus {\mathcal N}(L)$, see \cite[Theorem 8.20]{Go85}. Hence, because the constant functions are in ${\mathcal N}(L)$, the operator $(L, \D)$ is essentially m-dissipative in $\H$.

We verify assumptions $(H1)$-$(H4)$ as follows.

\paragraph{Proof of $(H1)$:}
Let $f\in \H$. Then $\pi_1 f\in L^2(\mu_1)$ with $\mu_1(\pi_1 f)=0.$ Let $\{g_n\}_{n\ge 0}\subset C_c^\infty(\R^{d_1})$ such that $\mu_1(g_n)=0$ and $\mu_1(|g_n-\pi_1 f|^2)\to 0.$ Let $\tt g_n(x,y) = g_n(x).$ Then $\tt g_n\in \D$, $\mu(|\tt g_n- \pi_1 f|^2)=\mu_1(|g_n-\pi_1 f|^2)\to 0$ and
$$\lim_{n,m\to\infty} \mu(|\tt g_n-\tt g_m|^2 +|S\tt g_n -S\tt g_m|^2) = \lim_{n,m\to\infty} \mu(|\tt g_n-\tt g_m|^2)=0.$$
Thus, $\{\tt g_n\}_{n\ge 1}$ is a Cauchy sequence in $\D(S)$ with $S g_n=0$, and converges to $\pi_1 f$ in $L^2(\mu)$. Therefore, $\pi_1 f\in \D(S)$ and $S\pi_1 f=0$  since the operator is closed.

\paragraph{Proof of $(H2)$:} For any $f\in \D$, we have $\pi_1 f\in \D$ depending only on the first component. So, $\pi_1 \D\subset \D\subset \D(A).$  Since $(\pi_1 f)(x,y)=\pi_1 f(x)$ only depends on $x$, by the definitions of $A$ and $\pi_1$, we have
\begin{align*}
-(\pi_1A\pi_1)f(x,y) & = \int_{\R^{d_2}} \<Q\nn^{(2)} V_2(y'), \nn^{(1)} \pi_1 f(x)\>\mu_2(\d y') = \<\mu_2(Q\nn^{(2)} V_2),\nn^{(1)} \pi_1f(x)\>=0,
\end{align*}
where the last step is due to
$V_2(y)=\Phi(|y|^2)$   and $|\nn V_2|\in L^1(\mu_2)$ according to $(H)$. Then $(H2)$ holds.

\paragraph{Proof of $(H3)$:} It suffices to prove \eqref{HH3}. For the first inequality, we only need to find out a bounded measurable function $K$ such that
\beq\label{5.2} SA\pi_1 f= K A\pi_1 f,\ \ f\in \D,\end{equation} since this implies
\beg{align*} BS& = (I+(A\pi_1)^*A\pi_1)^{-1} (A\pi_1)^* S = (I+(A\pi_1)^*A\pi_1)^{-1} (SA\pi_1)^* \\
&= (I+(A\pi_1)^*A\pi_1)^{-1} (KA\pi_1)^*  = BK,\end{align*} so that by
 $\|B\|\le 1$ we have
$$|\<BS\pi_2f, \pi_1f\> |=|\<BK\pi_2 f, \pi_1 f\>|\le \|K\|_\infty \|\pi_2 f\|\cdot\|\pi_1f\|.$$
Now for any $f\in \D$, \eqref{QP} implies
\beg{align*} &(S A\pi_1 f)(x,y) = S\<Q\nn^{(2)} V_2, \nn^{(1)} \pi_1 f\>(x,y)\\
&= (\DD^{(2)}-\<\nn^{(2)} V_2, \nn^{(2)}\cdot\>) \sum_{i=1}^{d_1} \big(2 \Phi'(|y|^2)(Qy)_i \pp_{x_i} \pi_1 f(x)\big)\\
&= 2 \sum_{i=1}^{d_1} \Big(\Phi''(|y|^2)(2d_2 - 4\Phi'(|y|^2) |y|^2 + 4) -2 \Phi'(|y|^2)^2 + 4 \Phi'''(|y|^2) |y|^2\Big) (Qy)_i \pp_{x_i} \pi_1 f(x)\\
&= 2 H(|y|^2) \<Q\nn^{(2)} V_2(y),\nn^{(1)} \pi_1 f(x)\> = 2 H(|y|^2) (A\pi_1 f)(x,y),\end{align*}
where
$$ H(r):= \ff{2r\Phi'''(r) +(d_2+2)\Phi''(r)}{\Phi'(r)}-\Phi'(r)-2r\Phi''(r),\ \ r>0,$$ is bounded according to $(H)$. Then \eqref{5.2} holds for
some bounded function $K$.

To prove the second inequality in \eqref{HH3}, we consider the operator $G:= -\pi_1A^2\pi_1= (A\pi_1)^*A\pi_1$ on $\D$.
By the definitions of $A$ and $\pi_1$, we have
\beq\label{5*3}   \beg{split} (G f)(x,y)=(Gf)(x)= \int_{\R^{d_2}}&-\Hess_{\pi_1 f}\Big(Q\nn^{(2)} V_2(y'), Q\nn^{(2)} V_2(y')\Big)(x)\\
&\Hess_{V_2} \Big(Q^*\nn^{(1)} V_1(x), Q^*\nn^{(1)} \pi_1 f(x)\Big)(y')\mu_2(\d y').\end{split}\end{equation} Then \eqref{QP} implies
\beg{align*} &\int_{\R^{d_2}}\Hess_{\pi_1 f}\Big(Q\nn^{(2)} V_2(y), Q\nn^{(2)} V_2(y)\Big)(x)\mu_2(\d y)\\
 &=4\sum_{i,j=1}^{d_1}  \int_{\R^{d_2}} \big(\pp_{x_i}\pp_{x_j} \pi_1 f\big)(x)\Phi'(|y|^2)^2 (Qy)_i(Qy)_j \mu_2(\d y)\\
&= 4\sum_{i,j=1}^{d_1}\sum_{k=1}^{d_2} \int_{\R^{d_2}} \big(\pp_{x_i}\pp_{x_j}  \pi_1 f\big)(x)\Phi'(|y|^2)^2Q_{ik}Q_{ik} y_k^2 \mu_2(\d y)\\
 &= \ff {4  } {d_2}\sum_{i,j=1}^{d_1} \int_{\R^{d_2}} (QQ^*)_{ij}\big(\pp_{x_i}\pp_{x_j}\pi_1 f\big)(x)\Phi'(|y|^2)^2 |y|^2 \mu_2(\d y) \\
& = \ff {\mu_2(|\nn V_2|^2)} {d_2}\sum_{i,j=1}^{d_1} (QQ^*)_{ij}   \big(\pp_{x_i}\pp_{x_j} \pi_1 f\Big)(x).\end{align*} Similarly,
\beg{align*} &\int_{\R^{d_2}} \Hess_{V_2} \Big(Q^*\nn^{(1)} V_1(x), Q^*\nn^{(1)} \pi_1 f(x)\Big)(y) \mu_2(\d y)\\
&= \<Q^*\nn^{(1)} V_1(x), Q^*\nn^{(1)} \pi_1 f(x)\> \int_{\R^{d_2}} 2\Phi'(|y|^2) + \ff{4\Phi''(|y|^2)|y|^2}{d_2}\mu_2(\d y)\\
&= \ff { \<Q^*\nn^{(1)} V_1(x), Q^*\nn^{(1)} \pi_1 f(x)\>} {d_2} \int_{\R^{d_2}} \DD^{(2)} V_2 (y) \mu_2(\d y)\\
& = \ff{\mu_2(|\nn^{(2)} V_2|^2)}{d_2}  \<Q^*\nn^{(1)} V_1(x), Q^*\nn^{(1)} \pi_1 f(x)\>.   \end{align*} Therefore, letting $N(V_2)=  \ff{\mu_2(|\nn^{(2)} V_2|^2)} {d_2}$ which is a positive constant according to $(H)$, we obtain
\beq\label{5*4}   (G f)(x,y)= (Gf)(x)=  -N(V_2) \sum_{i,j=1}^{d_1}(QQ^*)_{ij} \big\{\pp_{x_i}\pp_{x_j} - (\pp_{x_j} V_1)(x) \pp_{x_i}\big\}\pi_1 f(x).  \end{equation}
This enables us to provide the following assertion.

\beg{lem} \label{LL}   $(I+G)(\D)$ is dense in $\H$, so that $(-G,\D)$ is essentially m-dissipative  (equivalently, essentially self-adjoint) on $\H$.
\end{lem}

\begin{proof} First recall that for densely defined, symmetric and dissipative linear operators on a Hilbert space, the property of being essential m-dissipative is equivalent to essential self-adjointness. Consider the operator $(T,C_c^\infty(\R^{d_1}))$ on the Hilbert space $L^2(\mu_1)$ defined by
\beq\label{T} T:= \sum_{i,j=1}^{d_1}(QQ^*)_{ij} \big\{\pp_{x_i}\pp_{x_j} - (\pp_{x_j} V_1)(x) \pp_{x_i}\big\}.\end{equation}    Using integration by parts formula we have
\begin{align*}
\<Th, g \>_{L^2(\mu_1)} = - \mu_1(\<QQ^*\nabla^{(1)} h,  \nabla^{(1)} g\>),\ \ f\in C_c^\infty(\R^{d_1}), g\in C^\infty(\R^{d_1}).
\end{align*}   By \cite[Theorem 7]{BKR97} or \cite[Theorem 3.1]{Wie85} our assumptions in $(H)$ imply that $(T,C_c^\infty(\mathbb{R}^{d_1}))$ is essentially self-adjoint (hence, essentially m-dissipative) on $L^2(\mu_1)$. Therefore, $L^2(\mu_1) = \overline{{\mathcal R}(T)} \oplus {\mathcal N}(T)$. By \eqref{WPi} the null space ${\mathcal N}(T)$ consists of the constant functions only. Hence $(T,C_c^\infty(\mathbb{R}^{d_1}))$ restricted to $\H_1 = \{g\in L^2(\mu_1): \mu_1(g)=0\}$ is also essentially self-adjoint. Thus, $(I+G)(\D)$ is dense in $\H$, because $\H = \H_1 \oplus \H_2$ and $G$ acts trivial on $\H_2$.
%\begin{align} \label{LV}
%\<(I-G)f,g\>=0 \quad \mbox{for all } f \in \D.
%\end{align}
%It suffices  to show that $g=0$. Choose $\phi\in C_0^\infty(\R^d)$ with $0\le \phi\le 1$ and $\phi|_{B(0,1)}=1,$ and let
%$\phi_n(y)= \phi(y/n)$ for $n\ge 1$. For any function $f$ on $\R^{d_1+d_2}$ we define
%$$(f\star \phi_n)(x,y)= f(x,y)\phi_n(y),\ \ (x,y)\in\R^{d_1+d_2},$$ and we extend a function $f$ on $\R^{d_1}$ onto $\R^{d_1+d_2}$ by setting
% $\tt f(x,y) = f(x)$. Then from  \eqref{LV} and dominated convergence give that
%\begin{align*}
%&\<(I-N(V_2)\,T)f, g_1\>_{L^2(\mu_1)}=\<\tt f,g\>
%-  \<G\tt f,g\> \\
%&=\lim_{n\to\infty} \Big\{\<\tt f\star \phi_n,g\> -   \<G(\tt f\star \phi_n),  g \> \Big\}\\
%&=\lim_{n\to \infty} \<(I-G)(\tt f\star \phi_n),g\> = 0,\ \ f \in \pi_1\D,
%\end{align*}
% where $g_1=\pi_1 g$.
%Since $(I-N(V_2)\,T)(\pi_1\D)$ is dense in $\H_1:=\{f\in L^2(\mu_1): \mu_1(f)=0\},$  this implies  $g_1=0$.  Therefore,
%\begin{align*}
%\<Gf,g\>= N(V_2)\, \<T\pi_1 f, g_1\>_{L^2(\mu_1)}=0,\ \ f\in\D.
%\end{align*}
%Combining this with   \eqref{LV} we obtain  $\<f,g\>=0$ for all $f \in \D$. Hence $g=0$ as desired.
\end{proof}

Now we continue to prove the second inequality in \eqref{HH3}. Let $f\in \D$ and $g= (I+G)f$. As in \eqref{5*3}, by the definitions of $A$ and $\pi_1$  we have
\begin{align*}
&(A^2\pi_1f)(x,y) \\= & \,\Hess_{\pi_1 f}(Q\nn^{(2)} V_2(y'), Q\nn^{(2)} V_2(y))(x)
- \Hess_{V_2} (Q^*\nn^{(1)} V_1(x), Q^*\nn^{(1)} \pi_1 f(x))(y).
\end{align*}
So,
\beq \label{EL} \beg{split}
\| A^2 \pi_1 f  \| &\leq  \left\|  Q\nabla^{(2)} V_2   \right\|_{L^4(\mu_2)}^2 \left\|  (\nabla^{(1)})^2  \pi_1 f  \right\|_{L^2(\mu_1)}\\
\quad &+  \left\|  (\nabla^{(2)})^2 V_2   \right\|_{L^2(\mu_2)}   \left\|  |Q^*\nabla^{(1)}  V_1 |\cdot  |Q^*\nabla^{(1)}  \pi_1 f | \right\|_{L^2(\mu_1)}.
\end{split}\end{equation}
Due to    \eqref{5*4} and \eqref{T} we see that   $\pi_1 f$  solves the elliptic equation
\begin{align*}
\pi_1 f -  N(V_2)T\pi_1 f =\pi_1 g \quad \mbox{in } L^2(\mu_1).
\end{align*}
By applying the elliptic \textit{a priori} estimates   from \cite[(2.2) and  Lemma 8]{DMS13} (or see \cite[Section 5.1]{GS} for corresponding proofs including domain issues) to the right hand side of  \eqref{EL} we conclude
\beq\label{QQ}
\|(BA)^*g\|_{L^2(\mu)}    \leq c  \,\| \pi_1 g \|_{L^2(\mu_1)}  \leq c  \,\|g\|_{L^2(\mu)}
\end{equation}
for some  constant $c\in (0,\infty)$  only depending on  $V_1$ and $V_2$. According to Proposition \ref{P1.14} and Lemma \ref{LL}, this implies the second inequality in \eqref{HH3}. In conclusion, assumption $(H3)$ holds.

\paragraph{Proof of $(H4)$:} Let $f\in \D(L)$.
%For any $f\in \D(L)$, there exists a sequence $f_n\in \D$ such that $f_n\to f$ and $Lf_n\to Lf$ in $L^2(\mu)$.
Since $\mu(f)=0$, we have
$$\gg_1:={\rm ess}_\mu\inf f\le 0,\ \ \gg_2:={\rm ess}_\mu\sup f\ge 0.$$ Since $\D$ is a core of $(L,\D(L))$, we may take $\{g_n\}_{n\ge 1}\subset \D$ such that
$g_n\to f$ and $L g_n\to Lf$ in $L^2(\mu)$. To control $\|g_n\|_{\rm osc}$, for any $n\ge 1$ we take   $h_n\in C^\infty(\R)$ such that $0\le h_n'\le 1$ and
$$h_n(r)=\beg{cases} r\ &\text{for}\ r\in [\gg_1,\gg_2],\\
  \gg_1-\ff 1 {2n} \ &\text{for} \ r\le \gg_1 - \ff 1 n,\\
  \gg_2+\ff 1 {2n},\  &\text{for} \ r\le \gg_2 + \ff 1 n.\end{cases}$$
Then $f_n:= h_n(g_n)\to f$ in $L^2(\mu)$,
\beg{align*} &\limsup_{n\to \infty} \<-Lf_n,f_n\>= \limsup_{n\to \infty} \mu(h_n'(g_n)^2|\nn^{(2)}g_n|^2)\\
&\le \limsup_{n\to \infty} \mu(|\nn^{(2)}g_n|^2)=
\limsup_{n\to \infty} \<-Lg_n, g_n\>=\<-L f,f\>,\end{align*}  and
$$\limsup_{n\to \infty} \|f_n\|_{\rm osc}\le \limsup_{n\to \infty} \Big(\gg_2-\gg_1+\ff 1 n\Big)=\gg_2-\gg_1= \|f\|_{\rm osc}.$$
%Similarly,  for any $f\in\D(A\pi_1)$, we take $\{g_n\}_{n\ge 1}\subset\D$ such that $g_n\to f$ and $A\pi_1 g_n\to A\pi_1 f$ in $L^2(\mu)$. Since    $h_n'\in %C_b(\R)$ with $h_n'(f)=1$,   $h_n'(g_n)$ is bounded and converges to $1$ in $L^2(\mu)$. Thus, by dominated convergence, $f_n:=h_n(g_n)$ satisfies
%\beg{align*} &\limsup_{n\to\infty} \langle A\pi_1 f_n, A\pi_1 f_n\rangle\\ =  &\limsup_{n\to\infty}\int_{\R^{d_1+d_2}} \bigg(\int_{\R^{d_2}} h_n'(g_n)(x,y)\<Q\nn^{(2)}V_2(y'),\nn^{(1)}g_n(x,y)\>\mu_2(\d y)\bigg)^2\mu(\d x,\d y')\\
%= &\limsup_{n\to\infty}\int_{\R^{d_1+d_2}} \bigg(\int_{\R^{d_2}}  \<Q\nn^{(2)}V_2(y'),\nn^{(1)}g_n(x,y)\>\mu_2(\d y)\bigg)^2\mu(\d x,\d y')\\
%= &\limsup_{n\to\infty}\<A\pi_1g_n, A\pi_1g_n\> =\langle A\pi_1 f, A\pi_1 f \rangle. \end{align*}
Therefore, we have verified assumption $(H4)$.

\

\beg{proof}[Proof of Theorem \ref{T1.1}] It remains to prove \eqref{WC} for $\xi$ in \eqref{WP*}. Let $\Psi(f)= \|f\|_{\rm osc}^2$. The condition  \eqref{phi}  is obvious by the definition of $\pi_1$ and the $L^\infty(\mu)$-contraction of the Markov semigroups $P_t$ and $\e^{-t G}$.   Since we have verified assumptions $(H1)$-$(H4)$, by Theorem \ref{T2.1} it suffices to prove the weak Poincar\'e inequalities
\beq\label{WP1'} \|\pi_1 f\|^2\le c\aa_1(r) \|A\pi_1 f\|^2 + r\Psi(\pi_1 f),\ \ r>0, f\in \D(A\pi_1),\end{equation}
\beq\label{WP2'} \|\pi_2 f\|^2\le c\aa_2(r) \<Sf,f\>  + r\Psi(f),\ \ r>0, f\in \D \end{equation} for some constant $c\in (0,\infty).$

Recall that for any $f\in \D$  we have
$$(\pi_1f)(x,y)=  \int_{\R^{d_2}} f(x,y)\mu_2(\d y).$$
%Then $\pi_1 f\in C_c^\infty(\R^{d_1})$ and $\mu_1(\pi_1 f)=0$. So,   \eqref{WPi} for $i=1$ implies
%\beq\label{L1} \|\pi_1 f\|^2 =\mu_1(\pi_1 f^2) \le \aa_1(r) \mu_1(|\nn^{(1)} \pi_1 f|^2) + r \|\pi_1 f\|_{\rm osc}^2,\ \ r>0.\end{equation}
%Obviously, $\|\pi_1 f\|_{\rm osc}\le \|f\|_{\rm osc}$, and by 
By $V_2(y)=\Phi(|y|^2)$ we obtain 
\beg{align*} &\|A\pi_1 f\|^2 = \int_{\R^{d_1+d_2}} \<Q\nn^{(2)} V_2(y), \nn^{(1)} \pi_1 f(x)\>^2 \mu(\d x,\d y)\\
& = \ff 4 {Z(V_2)} \sum_{i,j=1}^{d_1} \int_{\R^{d_1}} (\pp_{x_i} \pi_1 f(x))(\pp_{x_j} \pi_1 f(x)) \mu_1(\d x) \int_{\R^{d_2}} \Phi'(|y|^2)^2(Q y)_i(Qy)_j \e^{-\Phi(|y|^2)}\d y\\
&= \ff{4}{Z(V_2)} \sum_{i,j=1}^{d_1}\sum_{k=1}^{d_2}Q_{ik}Q_{jk}  \int_{\R^{d_1}} (\pp_{x_i} \pi_1 f(x))(\pp_{x_j}\pi_1 f(x)) \mu_1(\d x) \int_{\R^{d_2}} \Phi'(|y|^2)y_k^2 \mu_2(\d y)\\
&=\ff{4 \int_{\R^d} |y|^2 \Phi'(|y|^2)^2 \mu_2(\d y)}{Z(V_2) d_2} \mu_1(|Q^*\nn^{(1)} \pi_1 f|^2).  \end{align*} Since $QQ^*$ is invertible, $0 < Z(V_2) < \infty$, and
$$0<4\int_{\R^{d_2}} |y|^2 \Phi'(|y|^2) \mu_2(\d y) = \mu_2(|\nn^{(2)}V_2|^2)<\infty$$ by $(H)$, this implies 
\beq\label{*M} \ff 1 c \mu_1(|\nn^{(1)} \pi_1 f|^2)\le \|A\pi_1 f\|^2\le c\mu_1(|\nn^{(1)}\pi_1f|^2),\ \ \ f\in\D \end{equation} for some constant $1< c < \infty$. 
So, $f\in \D(A\pi_1)$ implies that $\pi_1f\in  H^{1,2}(\mu_1)$, the completion of $C_c^\infty(\R^{d_1})$ under the Sobolev norm $\|g\|_{1,2}:=\ss{\mu_1(g^2+|\nn^{(1)}g|^2)}$.  Combining this with inequality \eqref{WPi} for $i=1$ which naturally extends to $f\in H^{1,2}(\mu_1)$, we prove    \eqref{WP1'}.

Next, for the above $f$ and  $x\in\R^d$,  we have $\hat f_x:= f(x,\cdot)-\pi_1 f(x)\in C_c^\infty(\R^d) $,   $\mu_2(\hat f_x)=0$ and $\|\hat f_x\|_{\rm osc}\le \|f\|_{\rm osc}.$ Then  \eqref{WPi} for $i=2$ implies
$$ \mu_2(|\hat f_x|^2) \le \aa_2(r) \mu_2(|\nn^{(2)} f(x,\cdot)|^2) + r\|  f \|_{\rm osc}^2,\ \ r>0.$$
Combining this with
\beg{align*} &\int_{\R^{d_1}} \mu_2(|\hat f_x|^2)\mu_1(\d x)= \|f-\pi_1f\|^2=\|\pi_2f\|^2,\\
&\int_{\R^{d_1}} \mu_2(|\nn^{(2)} f(x,\cdot)|^2) \mu_1(\d x) =\int_{\R^{d_1+d_2}} |\nn^{(2)} f|^2 (x,y) \mu(\d x,\d y) =-\<Lf,f\>,\end{align*}
we prove \eqref{WP2'} for $c=1$.
\end{proof}

To prove Example 1.1, we need the following lemma which is implied by the proof of \cite[Example 1.4]{RW01}.

\beg{lem}\label{WL} Let $\mu_V(\d x)= \e^{-V(x)}\d x$ be a probability measure on $\R^d$. Then the weak Poincar\'e inequality
\beq\label{WPV} \Var_{\mu_V}(f)\le r\aa_V(r) \mu_V(|\nn f|^2) + r \|f\|_{\rm osc}^2,\ \ r>0, f\in C_b^1(\R^d)\end{equation} holds for some decreasing $\aa_V: (0,\infty)\to [0,\infty).$ In particular: \beg{enumerate}
\item[$(1)$] If $V(x)\sim k |x|^\dd$ or $V(x)\sim k (1+|x|^2)^{\ff\dd 2}$ for some constants $k,\dd>0$, then \eqref{WPV} holds with
$$\aa_V(r)= c\big(\log(1+r^{-1})\big)^{\ff{4(1-\dd)^+}\dd}$$ for some constant $c>0.$
\item[$(2)$] If $V(x)\sim \ff{d+p}2\log(1+ |x|^2)$ for some constant $p>0$, then \eqref{WPV} holds with
$$\aa_V(r)= c r^{-\theta(p)} $$ for some constant $c>0.$
\item[$(3)$] If $V(x)\sim \ff d 2 \log(1+ |x|^2) + p \log\log (\e+ |x|^2)$ for some constant  $p>1$, then \eqref{WPV} holds with
$$\aa_V(r)= c_1\e^{c_2 r^{-\ff 1 {p-1}}}  $$ for some constant $c_1, c_2>0.$
\end{enumerate} \end{lem}

\beg{proof}[Proof of Example 1.1.] We only  consider case {\bf (A)} and the assertions in the other two cases can be verified in the same way.

By Lemma \ref{WL}, \eqref{WP1} holds for
\beq\label{AA1} \aa_1(r)= c\big(\log(\e + r^{-1})\big)^{\ff{4(1-\dd)^+}\dd}\end{equation}  for some constant $c>0$. Moreover, for case $(A_1)$, \eqref{WP2} holds for
$$\aa_2(r)= c'\big(\log(\e+r^{-1})\big)^{\ff{4(1-\vv)^+}{\vv}}.$$ Then for a constant $c_2>0$, there exists constants $\kk_1,\kk_2>0$ such that  the inequality
\beq\label{ETP}c_2t\ge  \aa_1(r)^2\aa_2\Big(\ff r{\aa_1(r)^2}\Big) \log\ff 1 r\end{equation}  implies
$$r\le \kk_1\exp\Big(-\kk_2t^{\ff{\dd\vv}{\dd\vv+ 8\vv(1-\dd)^++ 4\dd(1-\vv)^+}}\Big).$$ Therefore, the desired assertion  follows from \eqref{WP*}.

For case $(A_2)$ we may take
$$\aa_2(r)= c' r^{-\theta(p)}$$ for some constant $c'>0$. Then for a constant $c_2>0$, there exists constants $\kk>0$ such that  the inequality
\eqref{ETP}  implies
$$r\le \kk t^{-\ff 1 {\theta(p)}} \big(\log(\e+t)\big)^{\ff{8(\theta(p)+1)(1-\dd)^+ +\dd}{\theta(p)\dd}},$$  so that  the desired assertion  follows from \eqref{WP*}.

Finally, for case $(A_3)$ we may take
$$\aa_2(r)= c' \exp\Big(c'' r^{-\ff 1 {p-1}}\Big)$$ for some constants $c',c''>0$. Then for a constant $c_2>0$, there exists constants $\kk >0$ such that  the inequality \eqref{ETP} implies
$$r\le \kk  \big(\log(\e+t)\big)^{-(p-1)} \cdot \big(\log\log(\e^2+t)\big)^{\ff{8(1-\dd)^+}\dd},$$ so that the desired assertion   follows from \eqref{WP*}.
\end{proof}

\end{document}